\documentclass[12pt]{article}
\usepackage{amssymb,amsmath,latexsym,amsfonts,amsthm,mathrsfs}
\usepackage{enumerate}
\newtheorem{thm}{Theorem}
\newtheorem*{defn}{Definition}
\newtheorem{prop}{Proposition}
\newtheorem*{ex}{Example}
\newtheorem{exa}{Example}
\newtheorem*{claim}{Claim}
\newtheorem{claimm}{Claim}

\newtheorem*{ack}{Ackowledgements}
\newtheorem{lemma}{Lemma}

\newtheorem{cor}{Corollary}
\newtheorem*{note}{Notation}

\newcommand{\N}{\mathbb{N}}
\newcommand{\R}{\mathbb{R}}

\begin{document}

\author{J. G. Mijares}

\title{A notion of selective ultrafilter corresponding to topological Ramsey spaces.}

\date{}

\maketitle

\abstract{We introduce the relation of \textit{almost-reduction} in an arbitrary topological  Ramsey space $\mathcal{R}$ as a generalization of the relation of almost-inclusion on $\N^{[\infty]}$. This leads us to a type of ultrafilter $\mathcal{U}$ on the set of first approximations of the elements of $\mathcal{R}$ which corresponds to the well-known notion of \textit{selective ultrafilter} on $\N$. The relationship turns out to be rather exact in the sense that it permits us to lift several well-known facts about selective ultrafilters on $\mathbb{N}$ and the Ellentuck space  $\mathbb{N}^{[\infty]}$ to the ultrafilter $\mathcal{U}$ and the Ramsey space $\mathcal{R}$. For example, we prove that the Open Coloring Axiom holds on $L(\mathbb{R})[\mathcal{U}]$, extending therefore the result from \cite{DiT} which gives the same conclusion for the Ramsey space $\N^{[\infty]}$.}

\bigskip

\section{Preliminars.}

We follow \cite{Tod} in describing what a \textit{topological Ramsey space} is, rather than the earlier reference \cite{CaS} where a slightly different definition is given.  Consider triplets of the form $(\mathcal{R}, \leq, (r_{n})_{n\in\mathbb{N}} )$, where $\mathcal{R}$ is a set, $\leq$ is a pre-order on $\mathcal{R}$  and for every $n\in\mathbb{N}$, $r_{n}: \mathcal{R}\rightarrow \mathcal{AR}_{n}$ is a function with range $\mathcal{AR}_{n}$. If $A\leq B$ we say that $A$ is a \textit{reduction} of $B$; and for each $A\in \mathcal{R}$, we say that $r_{n}(A)$ is \textit{the} $n$th \textit{approximation of} $A$. We will assume that the following is satisfied:

\begin{itemize}
\item[{(A1)}]For any $A\in \mathcal{R}$, $r_{0}(A) = \emptyset$.
\item[{(A2)}]For any $A,B\in \mathcal{R}$, if $A\neq B$ then $(\exists n)r_{n}(A)\neq r_{n}(B)$.
\item[{(A3)}]If $r_{n}(A) = r_{m}(B)$ then $n = m$ and $(\forall i<n)r_{i}(A) = r_{i}(B)$.
\end{itemize}

These three assumptions allow us to identify each $A\in \mathcal{R}$ with the sequence $(r_{n}(A))_{n}$ of its approximations. In this way, if we consider the space $\mathcal{AR}=\bigcup_{n} \mathcal{AR}_{n}$ with the discrete topology, we can identify  $\mathcal{R}$ with a subspace of the (metric) space $\mathcal{AR}^{\mathbb{N}}$ (with the product topology) of all the sequences of elements of $\mathcal{AR}$. Via this identification, we will regard $\mathcal{R}$ as a subspace of $\mathcal{AR}^{\mathbb{N}}$, and we will say that $\mathcal{R}$ is {\it metrically closed} if it is a closed subspace of $\mathcal{AR}^{\mathbb{N}}$.

Also, for $a\in\mathcal{AR}$ we define the \textit{lenth} of $a$, $|a|$, as the unique $n$ such that $a = r_{n}(A)$ for some $A\in \mathcal{R}$.  

\medskip

We will further identify $a$ with the sequence $\{ r_{i}(A)\}_{i\leq n}$.  So, if $a = r_{n}(A)$ and $a' = r_i(A)$ for $i\leq n$ then we write $a' = r_i(a)$ (that is, we are extending the domain of the function $r_i$ to the set of $a\in\mathcal{AR}$ with $i\leq |a|$). In this case we also write $a'\sqsubseteq a$ and say that $a'$ is an initial segment of $a$.

\bigskip

We shall also consider on $\mathcal{R}$ the \textit{Ellentuck type neighborhoods} $$[a,A] = \{B\in \mathcal{R} : (\exists n)(a = r_{n}(B))\ \ \mbox{and} \ \ B\leq A\}$$
where $a\in\mathcal{AR}$ and $A\in \mathcal{R}$. If $[a,A]\neq\emptyset$ we will say that $a$ is \textit{compatible} with $A$ (or $A$ is compatible with $a$). Let $\mathcal{AR}(A) = \{a\in \mathcal{AR} : a\ \ \mbox{is compatible with}\ \ A\}$.

\vspace{.25 cm}

We write $[n,A]$ for $[r_{n}(A),A]$, and $Exp(\mathcal{R})$ for the family of all the neighborhoods $[n,A]$. This family generates the natural "exponential" topology on $\mathcal{R}$ which is finer than the product topology.

\bigskip

\begin{defn} 
A set $\mathcal{X}\subseteq \mathcal{R}$ is \textbf{Ramsey} if for every neighborhood $[a,A]\neq\emptyset$ there exists  $B\in [depth_A(a),A]$ such that $[a,B]\subseteq \mathcal{X}$ or $[a,B]\cap \mathcal{X} = \emptyset$. A set $\mathcal{X}\subseteq \mathcal{R}$ is \textbf{Ramsey null} if for every neighborhood $[a,A]\neq\emptyset$ there exists  $B\in [depth_A(a),A]$ such that $[a,B]\cap \mathcal{X} = \emptyset$.
\end{defn}

\bigskip

\begin{defn}
We say that $(\mathcal{R}, \leq, (r_{n})_{n})$ is a (topological) \textbf{Ramsey space} if every subset of $\mathcal{R}$ with the Baire property is Ramsey and every meager subset of $\mathcal{R}$ is Ramsey null.
\end{defn}

\vspace{.25 cm}

In \cite{Tod} it is shown that A1, A2 and A3, together with the following three assumptions are conditions of suficiency for a triplet $(\mathcal{R}, \leq, (r_{n})_{n})$, with $\mathcal{R}$ metrically closed, to be a Ramsey space. (This fact is called \textit{Abstract Ellentuck theorem} in \cite{Tod}):

\newpage 

\begin{flushleft}
(A4)(\textit{Finitization}) There is a pre-order $\leq_{fin}$ on $\mathcal{AR}$ such that:
\end{flushleft}

    \begin{itemize}
    \item[{(i)}]$A\leq B$ iff
    $\forall n\exists m \ \ r_{n}(A)\leq_{fin} r_{m}(B)$.
    \item[{(ii)}]$\{b\in \mathcal{AR} : b\leq_{fin} a\}$ is finite, for every
    $a\in \mathcal{AR}$.
    \end{itemize}

\bigskip

Given $a$ and $A$, we define the \textbf{depth of} $a$ \textbf{in} $A$, $depth_{A}(a)$, as

\medskip

$depth_A(a) =\left\{\begin{array}{ll} min\{n :a\leq_{fin}r_n(A)\}&\mbox{, if it exists.}\\ -1&\mbox{, otherwise.}\end{array}\right.$

\medskip

\begin{itemize}

\item[{(A5)}](\textit{Amalgamation}) Given compatible $a$ and $A$ with $depth_{A}(a)=n$, the following holds:
\begin{itemize}
\item[{(i)}] $\forall B\in [n,A]\ \ ([a,B]\neq\emptyset)$. 
\item[{(ii)}] $\forall B\in [a,A]\ \ \exists A'\in [n,A]\ \ ([a,A']\subseteq [a,B])$.
\end{itemize}
\item[{(A6)}](\textit{Pigeon Hole Principle}) Given compatible $a$ and $A$ with $depth_{A}(a) = n$, for each partition $\phi:\mathcal{AR}_{|a|+1}\rightarrow \{0,1\}$ there is $B\in [n,A]$ such that $\phi$ is constant in $r_{|a|+1}[a,B]$.

\end{itemize}

\bigskip

\begin{exa}[Ellentuck's Space, see \cite{Ell}]
The prototypical example of a topological Ramsey space is $\mathcal{R} = \N^{[\infty]}$, the set of infinite subsets of $\N$, with $\leq\ =\ \subseteq$ and $r_n(A)$ equal to the set formed by the first $n$ elements of $A$. In this case, $\mathcal{AR} = \N^{[<\infty]}$, the set of finite subsets of $\N$,  $a\leq_{fin} b$ if and only if $a = \emptyset$ or $a\subseteq b$ and $max(a) = max(b)$; and A6 reduces to the classical pigeon hole principle for infinites subsets of $\N$.
\end{exa}

\begin{exa}[Milliken's Space, see \cite{Milk}]
Let $FIN = \N^{[<\infty]}\setminus\{\emptyset\}$. Let $\mathcal{R} = FIN^{[\infty]}$, the set of infinite block sequences of elements of $FIN$, i.e., infinite subsets of $FIN$ of the form $A = \{a_0, a_1, a_2, \dots\}$ with $max(a_i) < min(a_{i+1})$. We write  $A\leq B$ if and only if $A\subset FU(B)$, where $FU(x)$ is the set of finite unions of elements of $x$, for any $x\subseteq FIN$. As in the previous example, let $r_n(A)$ be equal to the set formed by the first $n$ elements of $A$. In this case, $\mathcal{AR} = FIN^{[<\infty]}$, the set of finite block sequences of elements of $FIN$;  $a\leq_{fin} b$ if and only if $a = \emptyset$ or $a\subseteq FU(b)$ and $max(\bigcup a) = max(\bigcup b)$; and A6 reduces to Hindman's theorem \cite{Hin}.
\end{exa}

\section{The Relation of \textit{almost-reduction}.}

From now on we will assume that $(\mathcal{R}, \leq, (r_{n})_{n\in\mathbb{N}})$ satisfies A1-A6 and $\mathcal{R}$ is metrically closed; hence in virtue of the abstract Ellentuck theorem it is a topological Ramsey space.

\bigskip

In this section we introduce the following relation on $\mathcal{R}$, $$A\leq^*B\ \ \mbox{if and only if}\ \ \exists a\in\mathcal{AR}\ \ (\emptyset \neq [a,A]\subseteq [a,B]).$$

We call $\leq^*$ \textit{\textbf{relation of almost-reduction}} on $\mathcal{R}$. This is an abstraction of the relations of \textit{almost-inclusion} and  \textit{almost-condensation} (see \cite{Bla}), for elements of $\N^{[\infty]}$, $FIN^{[\infty]}$, respectively; it is also an abstraction of the relation of \textit{being almost a subspace} for elements of the topological Ramsey space  $FIN_k^{[\infty]}$ (see \cite{Tod} for a definition of this space). The name \textit{almost-reduction} is related to the fact that $A\leq B$ (i.e, $A$ is -completely- a \textit{reduction} of $B$) if and only if $[0,A]\subseteq [0,B]$. So roughly speaking, $A\leq^*B$ tells us that  $A$ is a reduction of $B$ ``from some point on''.

\vspace{.25 cm}

Our porpuse now is to show that  $\leq^*$ is a $\sigma$-closed pre-ordering on $\mathcal{R}$ (see theorem \ref{thmSigma} below), but first we need to understand $\leq^*$ in terms of the finite approximations to elements of $\mathcal{R}$.

\medskip

Given $a\in\mathcal{AR}(A)$ with $depth(a)\geq 0$, notice that 

$$\exists B\in[0,A]\ \ ([a,B]\neq\emptyset).$$

\begin{defn}
For $A\in\mathcal{R}$ and $a\in\mathcal{AR}(A)$ with $depth(a)\geq 0$, the {\bf depth 0} of $a$ in $A$ is defined as 
$$depth_A^0(a) = max\{n\leq depth_A(a): \forall b\leq_{fin}r_n(A)\ \exists B\in [b,A]\ \ ([a,B]\neq\emptyset)\}.$$
\end{defn} 

\bigskip

Notice that $depth(a) = 0$ if and only if $a = \emptyset$.

\begin{ex}
In Ellentuck's space, $depth_A^0(a) = n > 0$ if and only if the $n$-th element of $A$ is $min(a)$. In Milliken's space, if $A = \{x_1, x_2, \dots\}$ and $a\in\mathcal{AR}(A)$ then  $depth_A^0(a) = n > 0$ if and only if $x_n\subseteq \bigcup a$ and  $min(x_n) = min(\bigcup a)$; so $\forall i<n\ \ max(x_i) < min(\bigcup a)$. In both cases,  $depth_A^0(a)$ gives us a measure of which is \textbf{the ``least'' element of $A$ used to ``build'' $a$}. This is the intended idea in the general case.  
\end{ex}

\vspace{2 cm}

\textbf{Notice} that $depth_A^0(a)$ satisfies the following:

\begin{enumerate}
\item{$depth_A(a)\geq 0\ \rightarrow\ 0 \leq depth_A^0(a)\leq depth_A(a)$.}
\item{$depth_A(a)>0\ \rightarrow\ \ depth_A^0(a) = depth_A^0(r_1(a)) \geq 1$.}
\end{enumerate}

\bigskip

\begin {defn}
For $A\in\mathcal{R}$ and $a,b\in\mathcal{AR}(A)$, we write $a <_A b$ to mean $depth_A(a) < depth_A^0(b)$.

\bigskip

Also, for $A\in\mathcal{R}$ and $a\in\mathcal{A}$, $$\mathcal{AR}(A)/a = \{b\in\mathcal{AR}(A) : a <_A b\}.$$ We write $\mathcal{AR}(A)/n$ to denote $\mathcal{AR}(A)/r_n(A)$. And for every $i>0$ we write $\mathcal{AR}_i(A)/a$ to denote $ \mathcal{AR}_i(A)\cap(\mathcal{AR}(A)/a)$.

\bigskip

Also, for $a,b\in\mathcal{AR}(A)$, $$b/a = \{c\leq_{fin}b :  a <_A c \}.$$
\end{defn}

\medskip

\begin{lemma}\label{straight}
The following facts are straightforward:

\begin{enumerate}
\item{$A\leq B$ implies $A\leq^* B$.}
\item{$A\leq^*B\ \ \mbox{if and only if}\ \ \exists n\ \mathcal{AR}(A)/n\subseteq \mathcal{AR}(B)$.}

\item{$A\leq^* B$ if and only if $\exists n\ \mathcal{AR}_1(A)/n\subseteq\mathcal{AR}_1(B).$}
\end{enumerate}
\end{lemma}
\qed

\bigskip

\begin{thm}\label{thmSigma}
$(\mathcal{R},\leq^*)$ is a $\sigma$-closed  pre-ordered set.
\end{thm}
\begin{proof}
Clearly, $\leq^*$ is reflexive. So let us see that it is transitive. Suppose $A\leq^*B$ and $B\leq^*C$. So by lemma \ref{straight} there exist $n,m\in\omega$ such that 
\begin{enumerate}
\item[{(1)}]$\mathcal{AR}(A)/n\subseteq \mathcal{AR}(B)$.
\item[{(2)}]$\mathcal{AR}(B)/m\subseteq \mathcal{AR}(C)$.
\end{enumerate}

 Choose an integer $k\geq n$ as follows. If there exists $c\in\mathcal{AR}(A)/n$  such that $depth_B^0(c)\leq m$ then let $$k =  max\{depth_A(c) : c\in\mathcal{AR}(A)/n \ \ \&\ \ depth_B^0(c)\leq m\};$$ otherwise, let $k = n$.

\medskip

By (2), $k = n$ implies $\mathcal{AR}(A)/n\subseteq \mathcal{AR}(C)$ and we are done. And if  $k>n$, then  by (1) $\mathcal{AR}(A)/k\subseteq\mathcal{AR}(A)/n\subseteq \mathcal{AR}(B)$. Take $c\in\mathcal{AR}(A)/k$. Then $depth_A^0(c) > k$, and therefore $depth_A(c) > k$. By the choice of $k$, this implies $depth_B^0(c) > m$. So again, by (2), this proves $\mathcal{AR}(A)/k\subseteq \mathcal{AR}(C)$. In any case, $A\leq^* C$. 

\vspace{.25 cm}

Now to prove the $\sigma$-closedness, let $(A_n)_{n\geq 1}$ be a decreasing sequence in  $(\mathcal{R},\leq^*)$. Notice that there exists a sequence $(a_n)_{n\geq 1}$ satisfying 

\begin{enumerate}
\item[{i.}]$a_1\in\mathcal{AR}_1(A_1)$,  
\item[{ii.}]$(\forall n>1)\ \ a_n\in r_{|a_{n-1}|+1}[a_{n-1},A_1]$ and
\item[{iii.}]$(\forall n>1)\ \ a_n/a_{n-1}\subset\mathcal{AR}(A_n)$.
\end{enumerate}

\medskip

 Let  $B\in\mathcal{R}$ be such that for every $n\geq 1$, $r_n(B) = a_n$. Obviously $B\leq A_1$, and for  $n > 0$ we have $\mathcal{AR}(B)/n\subseteq\mathcal{AR}(A_{n+1})$ and so $(\forall n)\ B\leq^* A_n$.
\end{proof}

\bigskip

\begin{defn}

We say that  $\mathcal{D}\subseteq\mathcal{R}$ is {\em dense open} in  $(\mathcal{R},\leq^*)$ if:
\begin{enumerate}
\item[{i.}]for every $A\in\mathcal{R}$ there is $B\in\mathcal{D}$ such that $B\leq^* A$, and 
\item[{ii.}]for every $A, B\in\mathcal{R}$, $B\leq^* A$ and $A\in\mathcal{D}$ imply $B\in\mathcal{D}$.
\end{enumerate}

\end{defn}

\begin{cor}
For every sequence $(\mathcal{D}_n)_n$ of dense open subsets of $\mathcal{R}$, the set $\bigcap_n\mathcal{D}_n$ is dense in $(\mathcal{R},\leq^*)$.
\end{cor}
\begin{proof}
Fix $B\in\mathcal{R}$. Pick $A_0\in\mathcal{D}_0$ such that $A_0\leq^* B$; and for each $n>0$ take $A_n\in\mathcal{D}_n$ in such a way that $A_n\leq^* A_{n-1}$. Apply theorem \ref{thmSigma} to obtain $A\in\mathcal{R}$ such that $A\leq^* A_n$, for every $n$. Hence by transitivity $A\leq^* B$ and since each $\mathcal{D}_n$ is open, we have $A\in\bigcap_n\mathcal{D}_n$.
\end{proof}

\vspace{.25 cm}

We finish this section by proving the following interesting consequence of theorem \ref{thmSigma} and A6. This is an abstract version of two-dimensional Ramsey's theorem \cite{Ram}.

\begin{cor}\label{corRamsey}
 For every partition of $\mathcal{AR}_2$ into  two classes and for every $A\in\mathcal{R}$, there exists $B\leq A$ such that $\mathcal{AR}_2(B)$ lies in one single partition class.
\end{cor}
\begin{proof}
Let $A\in\mathcal{R}$ be given and consider an arbitrary partition $$\mathcal{AR}_2 = \mathcal{C}_0\cup\mathcal{C}_1.$$ For every $n\in\omega$ and every $a\in\mathcal{AR}_1(A)$ with $depth_A(a)=n$ define $$\mathcal{D}^a_n = \{B\in\mathcal{R} : B\leq A\ \ \mbox{and}\ \ (r_2[a,B]\subseteq\mathcal{C}_0\ \ \mbox{or}\ \ (r_2[a,B]\subseteq\mathcal{C}_1)\}.$$ Take $C\in [0,A]$. Notice that A6 (together with A5 in the case $[a,C] = \emptyset$) implies that there exists $B\in\mathcal{D}^a_n$ such that $B\leq^* C$. Then each $\mathcal{D}^a_n$ is dense open below $A$ and so for each $n$, $$\mathcal{D}_n = \bigcap_{depth_A(a)=n}\mathcal{D}^a_n$$ is also dense open below $A$ (by A4, the set $\{a\in\mathcal{AR}_1(A) : depth_A(a) = n\}$ is finite). Hence, by the previous corolary there exists $\bar{A}\in\bigcap_n\mathcal{D}_n$ such that $\bar{A}\leq^* A$. We can assume $\bar{A}\leq A$ without a loss of generality, because of the openness of $\bigcap_n\mathcal{D}_n$ (i.e., by the definition of $\leq^*$ there exists $b\in\mathcal{AR}$ such that $\emptyset\neq[b,\bar{A}]\subseteq [b,A]$, so we can choose any element of $[b,\bar{A}]$ instead of $\bar{A}$). By the definition of $\mathcal{D}_n$, for every $a\in\mathcal{AR}_1(\bar{A})$ there exist $i_a\in\{0,1\}$ such that $b\in\mathcal{C}_{i_a}$ for every $b\in r_2[a,\bar{A}]$. Consider the partition $c : \mathcal{AR}_1 \rightarrow \{0,1\}$ given by 

$$c(a) = i_a\mbox{,}\ \mbox{if}\ \ a\in\mathcal{AR}_1(\bar{A}).$$
Apply A6 to obtain $B\in [0,\bar{A}]$ such that $c$ is constant in $r_1[0,B] = \mathcal{AR}_1(B)$. So $B\leq A$ and there exists $i\in\{0,1\}$ such that $i_a = i$ for every $a\in\mathcal{AR}_1(B)$. This means $\mathcal{AR}_2(B)\subseteq\mathcal{C}_i$.
\end{proof}

\medskip

Corolary \ref{corRamsey}, together with A6, indicates that every topological Ramsey space behaves as  a \textit{happy family} or \textit{selective coideal} (see \cite{Mat} or \cite{Far}) .

\section{The forcing notion of \textit{almost-reduction}  and the corresponding ultrafilter.}

\bigskip

In this section we consider $\mathbb{A}  = (\mathcal{R},\leq^*)$ as a forcing notion. The following definition is related to the forcing-like nature of $\mathbb{A}$: 

\begin{defn}
A \textbf{first-approximation ultrafilter} is an ultrafilter $\mathcal{U}$ on $\mathcal{AR}_1$ which is generated by sets of the form $\mathcal{AR}_1(A)$ with $A\in\mathcal{R}$.
\end{defn}

\bigskip

We will assume that every first-approximation ultrafilter  $\mathcal{U}$ is {\em closed under finite changes}, i.e., if $\mathcal{AR}_1(A)\in\mathcal{U}$ then $\mathcal{AR}_1(A)/a\in\mathcal{U}$ for every $a\in\mathcal{AR}_1(A)$.

\bigskip

\begin{note}
For a family $\mathcal{U}$ of subsets of $\mathcal{AR}_1$, let $$\mathcal{R}_{\mathcal{U}} = \{A\in\mathcal{R} : \mathcal{AR}_1(A)\in\mathcal{U}\}.$$
\end{note}

\vspace{1 cm}

It turns out that if $G$ is the $\mathbb{A}$-generic filter then $$\mathcal{U}_G = \{X\subseteq\mathcal{AR}_1 : \exists A\in G\ \mathcal{AR}_1(A)\subseteq X\}$$ is a first-approximation ultrafilter which satisfies the following  very interesting property:

\bigskip

\begin{defn}
Let $\mathcal{U}$ be a first-approximation ultrafilter. We say that  $\mathcal{U}$ is \textbf{Ramsey} if for every partition of $\mathcal{AR}_2$ into  two classes and for every $A\in\mathcal{R}_{\mathcal{U}}$, there exists $B\in\mathcal{R}_{\mathcal{U}}$ such that $B\leq A$ and $\mathcal{AR}_2(B)$ lies in one single partition class.
\end{defn}

\vspace{.25 cm}

So $\mathcal{U}$ is Ramsey if an abstract version of two-dimensional Ramsey's theorem, ``modulo $\mathcal{U}$'', holds (compare with corolary \ref{corRamsey} above). The following lemma  summarizes the main features of the forcing notion  $\mathbb{A}$. The argument is similar to the proof of the same fact for $(\N^{[\infty]},\subseteq^*)$.

\begin{lemma}
Forcing with $\mathbb{A}$ adjoins no new reals and if $G$ is the  $\mathbb{A}$-generic filter over some ground model $M$, then $\mathcal{U}_G$ as defined above is a Ramsey first-approximation ultrafilter in  $M[G]$.
\end{lemma}
\begin{proof} Since $\mathbb{A}$ is $\sigma$-closed (by theorem \ref{thmSigma}) then it adds no new reals. Also, notice that  $G = \mathcal{R}_{\mathcal{U}_G}$ in $M[G]$. Take $A\in G$ and also take, in $M[G]$, a partition $\phi : \mathcal{AR}_2 \rightarrow \{0,1\}$. All these objects are actually in the ground model so apply corolary \ref{corRamsey} to obtain $B\leq^* A$ such that $\phi$ is constant in $\mathcal{AR}_2(B)$. This means that the set $$\mathcal{D} = \{C\in\mathcal{R} : \phi\ \ \mbox{ is constant in}\ \  \mathcal{AR}_2(C)\}$$ is dense and hence $G\cap\mathcal{D} \neq \emptyset$.
\end{proof}

\bigskip

As pointed out in the proof of the previous lemma, for every ground model $M$ containing $\mathbb{A}$, $G = \mathcal{R}_{\mathcal{U}_G}$ in $M[G]$. Hence $M[G] = M[\mathcal{U}_G]$ . So we will consider $\mathcal{U}_G$ as \textbf{the} $\mathbb{A}$-\textbf{generic first-approximation ultrafilter over} $M$. And given $\mathcal{U}\subseteq\mathcal{AR}_1$ (in some ground  model $M$ containing $\mathbb{A}$), whenever we say that $\mathcal{U}$ is the $\mathbb{A}$-generic first-approximation ultrafilter over $M$ we mean that $\mathcal{U} = \mathcal{U}_G$ in $M[G]$ (where $G$ is the $\mathbb{A}$-generic filter over $M$).

\bigskip

\begin{defn}

Given a family  $\mathcal{U}$ of subsets of $\mathcal{AR}_1$, we say that  $\mathcal{D}\subseteq\mathcal{U}$ is {\em dense open} in $\mathcal{U}$ if $\mathcal{R}_{\mathcal{D}}$ is dense open in $(\mathcal{R}_{\mathcal{U}},\leq^*)$, i.e:
\begin{enumerate}
\item[{i.}]for every $A\in\mathcal{R}_{\mathcal{U}}$ there is $B\in\mathcal{R}_{\mathcal{D}}$ such that $B\leq^* A$, and 
\item[{ii.}]for every $A, B\in\mathcal{R}_{\mathcal{U}}$, $B\leq^* A$ and $A\in\mathcal{R}_{\mathcal{D}}$ imply $B\in\mathcal{R}_{\mathcal{D}}$.
\end{enumerate}

\end{defn}

\bigskip

\begin{defn}
Let $\mathcal{U}$ be a first-approximation ultrafilter. We say that  $\mathcal{U}$ is \textbf{selective} if for every sequence $(A_n)_n\subseteq\mathcal{R}_{\mathcal{U}}$ with $A_{n+1}\leq A_n$, there exists $B\in\mathcal{R}_{\mathcal{U}}$ such that $B\leq A_n$ for every $n$. So $B$ is a ``diagonalization'' of the sequence $(A_n)_n$ in $\mathcal{R}_{\mathcal{U}}$.
\end{defn}

\medskip

\begin{lemma}\label{lemmaEquiv}
Let $\mathcal{U}$ be a first-approximation ultrafilter. The following properties are equivalent:
\begin{enumerate}
\item[{(1)}]$\mathcal{U}$ is selective.
\item[{(2)}]For every sequence $(\mathcal{D}_n)_n$ of dense open subsets of $\mathcal{U}$, the set $\bigcap_n\mathcal{D}_n$ is dense in $\mathcal{U}$.
\item[{(3)}]$\mathcal{U}$ is Ramsey.
\end{enumerate}
\end{lemma}
\begin{proof} We are going to show $(1)\iff (2)$ and  $(2)\iff (3)$.

\vspace{1 cm}

$(1)\iff (2)$

\vspace{1 cm}

Let $\mathcal{U}$ be a selective first-approximation ultrafilter and for each $n\in\omega$ let $\mathcal{D}_n$ be a dense open subset of $\mathcal{U}$. Given $A\in\mathcal{R}_{\mathcal{U}}$, choose $A_n\in\mathcal{R}_{\mathcal{D}_n}$ such that $A_0\leq^* A$ and $A_{n+1}\leq^* A_n$ for each $n$. Actually, we can asume $A_{n+1}\leq A_n$ because $\mathcal{U}$ is closed under finite changes and each $\mathcal{D}_n$ is open. So by selectivity there is $B\in\mathcal{R}_{\mathcal{U}}$ such that $B\leq A_n$ for every $n$. By transitivity of $\leq^*$ we have $B\leq^* A$, and by openness of each $\mathcal{D}_n$ we also have $B\in\mathcal{R}_{\bigcap_n\mathcal{D}_n}$.

\vspace{1 cm}

Conversely, let $(A_n)_n\subseteq\mathcal{R}_{\mathcal{U}}$ with $A_{n+1}\leq A_n$ be given. Define $$\mathcal{D}_n = \{\mathcal{AR}_1(B)\in\mathcal{U} : B\leq A_n\}.$$ 

Each $\mathcal{D}_n$ is dense open in $\mathcal{U}$ and then $\bigcap_n\mathcal{D}_n$ is dense and therefore nonempty. For any $B\in\mathcal{R}_{\bigcap_n\mathcal{D}_n}$ we have $(\forall n)\ \ B\leq A_n$.

\vspace{2 cm}

$(2)\iff (3)$

\vspace{1 cm}

We proceed as in the proof of corolary \ref{corRamsey}. Let $A\in\mathcal{R}_{\mathcal{U}}$ be given and consider an arbitrary partition $$\mathcal{AR}_2 = \mathcal{C}_0\cup\mathcal{C}_1.$$ For every $n\in\omega$ and every $a\in\mathcal{AR}_1(A)$ with $depth_A(a)=n$ define $$\mathcal{D}^a_n = \{\mathcal{AR}_1(B)\in\mathcal{U} : B\leq A\ \ \mbox{and}\ \ (r_2[a,B]\subseteq\mathcal{C}_0\ \ \mbox{or}\ \ (r_2[a,B]\subseteq\mathcal{C}_1)\}.$$ Each $\mathcal{D}^a_n$ is dense open in $\mathcal{U}$ and so for each $n$, $$\mathcal{D}_n = \bigcap_{depth_A(a)=n}\mathcal{D}^a_n$$ is a dense open in  $\mathcal{U}$. So by (2) there exists $\bar{A}\in\mathcal{R}(\bigcap_n\mathcal{D}_n)$ such that $\bar{A}\leq^* A$. We can assume $\bar{A}\leq A$ without a loss of generality, because of the openness of $\bigcap_n\mathcal{D}_n$ and our assumption that $\mathcal{U}$ is closed under finite changes. Notice that by the definition of $\mathcal{D}_n$, for every $a\in\mathcal{AR}_1(\bar{A})$ there exist $i_a\in\{0,1\}$ such that $b\in\mathcal{C}_{i_a}$ for every $b\in r_2[a,\bar{A}]$. Consider the set $$X_0 = \{a\in\mathcal{AR}_1(\bar{A}) : i_a = 0\}.$$ Since $\mathcal{U}$ is an ultrafilter, one of the sets $X_0$ or $\mathcal{AR}_1(\bar{A})\setminus X_0$ is an element of $\mathcal{U}$ and then one of them contains $\mathcal{AR}_1(B)$ for some $B\in\mathcal{R}_{\mathcal{U}}$. So $B\leq A$ and by the definition of $X_0$ (or $\mathcal{AR}_1(\bar{A})\setminus X_0$, in any case) $\mathcal{AR}_2(B)$ is included in a single partition class.

\vspace{1 cm}

Conversely, let $(\mathcal{D}_n)_n$ be a sequence of dense open subsets of $\mathcal{U}$. Take $A_n\in\mathcal{D}_n$ such that $A_{n+1} \leq^* A_n$.  Define $c : \mathcal{AR}_2 \rightarrow \{0,1\}$ such that $$c(b) = 1\ \ \mbox{iff}\ \ r_1(b)\in\mathcal{AR}_1(A_0)\ \&\ \ b/r_1(b)\subset\mathcal{AR}(A_n)\ \ \forall n\leq depth_{A_0}(r_1(b)).$$ By Ramseyness, there exists $B\leq A_0$ in $\mathcal{R}_{\mathcal{U}}$ such that $c$ is constant on $\mathcal{AR}_2(B)$. Take any $a\in\mathcal{AR}_1(B)$. Notice that $\mathcal{AR}_1(B)\cap\bigcap\{\mathcal{AR}_1(A_n) : n\leq depth_{A_0}(a)\}$ belongs to $\mathcal{U}$. So there is $\bar{B}\in\mathcal{R}_{\mathcal{U}}$ such that $\bar{B}\leq B$ and $\bar{B}\leq A_n$ for every $n\leq depth_{A_0}(a)$. We can assume $[a,\bar{B}]\neq\emptyset$ since $\mathcal{U}$ is closed under finite changes. Take $b\in r_2 [a,\bar{B}]$. Then $c(b) = 1$, and therefore $c$ takes the constant value $1$ on $\mathcal{AR}_2(B)$. This means that for every $a\in\mathcal{AR}_1(B)$ we have $\mathcal{AR}(B)/a\subseteq\mathcal{AR}(A_n)$ for every $n\leq depth_{A_0}(a)$. This implies $(\forall n\in\omega)\ \ B\leq^* A_n$. And since each $\mathcal{D}_n$ is open we have $B\in\mathcal{R}_{\bigcap_n\mathcal{D}_n}$.

\end{proof}

\bigskip

We now show that the infinite-dimensional Ramsey-theoretic properties of selective ultrafilters on $\N$, discovered by Mathias \cite{Mat}, easily lift to the corresponding properties of first-approximation ultrafilters.

\bigskip

\begin{defn}
Let $\mathcal{U}$ be a first-approximation ultrafilter. We say that $\mathcal{X}\subseteq\mathcal{R}$ is $\mathcal{U}$-\textbf{Ramsey} if for every nonempty $[a,A]$ with $A\in\mathcal{R}_{\mathcal{U}}$ there exists $B\in [depth_A(a),A]\cap\mathcal{R}_{\mathcal{U}}$ such that $[a,B]\subseteq\mathcal{X}$ or $[a,B]\subseteq\mathcal{X}^{c}$.
\end{defn}

\bigskip

\begin{note}
For $a\in\mathcal{AR}$ and $B\in\mathcal{R}$ write $[a,B]$ \textbf{decides}  $\mathcal{X}$ if $$[a,B]\subseteq\mathcal{X}\ \mbox{or}\ \ [a,B]\subseteq\mathcal{X}^{c}.$$
\end{note}

\bigskip

\begin{lemma}\label{lemmaOpen}
Let $\mathcal{U}$ be a selective first-approximation ultrafilter. Then  every metric open subset of $\mathcal{R}$ is $\mathcal{U}$-Ramsey.
\end{lemma}

\begin{proof}
We proceed similarly to the proof of proposition 1.5 in \cite{Mat}. Let $\mathcal{X}\subseteq\mathcal{R}$ be open and fix a nonempty $[a,A]$ with $A\in\mathcal{R}_{\mathcal{U}}$. We will assume $a = \emptyset$ without a loss of generality. For every $b\in\mathcal{AR}$ define $$\mathcal{D}_b = \{\mathcal{AR}_1(B) : [b,B]\ \mbox{decides}\ \mathcal{X}\ \mbox{or}\ \forall C\in [0,B]\cap\mathcal{R}_{\mathcal{U}}\ \ ([b,C]\ \mbox{does not decide}\ \mathcal{X})\}.$$ 

Each $\mathcal{D}_b$ is dense open. So by lemma \ref{lemmaEquiv}, there exists $B\in\mathcal{R}_{\bigcap_b\mathcal{D}_b}$ such that $B\leq^* A$. We can assume $B\leq A$ without a loss of generality. Define $h : \mathcal{AR} \rightarrow\{0,1,2\}$ such that 

\begin{enumerate}
\item{$h(b) = 0$ iff $[b,B]\subseteq\mathcal{X}$.}
\item{$h(b) = 1$ iff $[b,B]\subseteq\mathcal{X}^{c}$.}
\item{$h(b) = 2$, otherwise.}
\end{enumerate}

\begin{claimm}\label{claim1}
If $h(b) = 2$ then there exists $C_b\in\mathcal{R}_{\mathcal{U}}\cap [b,B]$ such that $h(c) = 2$ for every $c\in r_{|b|+1}[b,C_b]$.
\end{claimm}
\begin{proof}[Proof of the claim \ref{claim1}.]

Define $\hat{h} : \mathcal{AR}_{|b|+1} \rightarrow\{0,1,2\}$ as $\hat{h}(c) = h(c)$. Then the result follows from A6 modulo $\mathcal{U}$ and the fact that $h(b) =2$.
\end{proof}

\begin{claimm}\label{claim2}
If $h(\emptyset) = 2$ then there exists $C\leq B$ in $\mathcal{R}_{\mathcal{U}}$ such that $\forall b\ \ ([b,C]\neq\emptyset \rightarrow h(b) = 2)$.
\end{claimm}

\begin{proof}[Proof of claim \ref{claim2}.]
For every $b\in\mathcal{AR}$ let $C_b$ be as in the previous lemma if $h(b) = 2$ and $[b,B]\neq\emptyset$; and $C_b = B$, otherwise. Let $C\in\mathcal{R}_{\mathcal{U}}$ be a diagonalization of $\{C_b\}_b$. Suppose there is $b$ such that $[b,C]\neq\emptyset$ and $h(b) \neq 2$. Assume $b$ is of minimal lenth. Notice that $b\neq\emptyset$ since we are assuming $h(\emptyset) = 2$, so let $b' = r_{|b|-1}(b)$. By minimality, we have $h(b') = 2$. But $b\in r_{|b'|+1}[b',C]\subseteq r_{|b'|+1}[b',C_{b'}]$ and hence by the previous claim $h(b) = 2$, contradicting our asumption on $b$. 
\end{proof}

\begin{claimm}\label{claim3}
$h(\emptyset) < 2$.
\end{claimm}
\begin{proof}[Proof of claim \ref{claim3}.]
Otherwise, let $C\leq B$ be as in the previous claim. Notice that $[0,C]$ does not decide $\mathcal{X}$, so there is $\bar{C}\leq C$ such that $\bar{C}\in\mathcal{X}$. But since $\mathcal{X}$ is a metric open there exists $b\sqsubset\bar{C}$ such that $[b] = \{D\in\mathcal{R} : b\sqsubset D\}\subseteq\mathcal{X}$. In particular, $[b,B]\subseteq\mathcal{X}$. But this means that $h(b) = 0$, contradicting our assumption on $h(\emptyset)$ in virtue of the previous claim.
\end{proof}

Claim \ref{claim3} implies  $[0,B]\subseteq\mathcal{X}$ or $[0,B]\subseteq\mathcal{X}^{c}$. This completes the proof of the lemma.
\end{proof}

\section{Other forcing notions related to topological Ramsey spaces.}\label{sec4}

In this section we use the methods of \cite{Mat} to analyse other abstract forcing notions related to topological Ramsey spaces. The results of this section will be of much help in section \ref{sect6} below, where we present an application of the previous results.
 
\bigskip

\subsection{Abstract Mathias forcing $\mathbb{M}$.}

\medskip

Let  $\mathbb{M}$ be the set of all the pairs $(a,A)$ such that $a\in\mathcal{AR}$ and $A\in\mathcal{R}$ and $[a,A]\neq\emptyset$. Order $\mathbb{M}$ as follows $$(a,A)\leq (b,B)\ \ \mbox{iff}\ \ [a,A]\subseteq [b,B].$$

\subsection{$\mathbb{P}_{\mathcal{U}}$}

Given a first-approximation ultrafilter $\mathcal{U}$ let $\mathbb{P}_{\mathcal{U}}$ be the set of all the pairs $(a,A)$ such that $a\in\mathcal{AR}$, $A\in\mathcal{R}_{\mathcal{U}}$ and $[a,A]\neq\emptyset$. Order $\mathbb{P}_{\mathcal{U}}$ as in the case of $\mathbb{M}$.

\bigskip

If $M$ is a transitive model of $ZF + DCR$ (for instance), we say that $g\in\mathcal{R}$ is $\mathbb{M}$-\textbf{generic} (resp. $\mathbb{P}_{\breve{\mathcal{U}}}$-\textbf{generic}) over $M$, if for every dense open subset $\mathcal{D}\in M$, of $\mathbb{M}$ (resp. $\mathbb{P}_{\breve{\mathcal{U}}}$), there exists a condition $(a,A)\in\mathcal{D}$ such that $g\in [a,A]$.

\bigskip

\begin{defn}
Let $\mathcal{U}$ be a selective first-approximation ultrafilter, $\mathcal{D}$ a dense open subset of  $\mathbb{P}_{\mathcal{U}}$ and $a\in\mathcal{AR}$. We say that $A$ {\em captures} $(a,\mathcal{D})$ if 

\begin{center}
$A\in\mathcal{R}_{\mathcal{U}}$, $[a,A]\neq\emptyset$ and $\forall B\in [a,A]\ \exists m > |a|\ \ (r_m(B),A)\in\mathcal{D}$.
\end{center} 
\end{defn}

\bigskip

\begin{lemma}
Let $\mathcal{U}$ be a selective first-approximation ultrafilter and $\mathcal{D}$ a dense open subset of $\mathbb{P}_{\mathcal{U}}$. Then for every $a\in\mathcal{AR}$ there exists $A$ which captures $(a,\mathcal{D})$.
\end{lemma}
\begin{proof}
Take any $B\in\mathcal{R}_{\mathcal{U}}$ with $[a,B]\neq\emptyset$. For every $b\in\mathcal{AR}(B)$ with $a\sqsubseteq b$ pick $C_b\in\mathcal{R}_{\mathcal{U}}$ such that $(b,C_b)\in\mathcal{D}$ if possible; otherwise set $C_b = B$. Let $C\in\mathcal{R}_{\mathcal{U}}\cap[a,B]$ be a diagonalization of $\{C_b\}_b$. Then for every $b$ as above such that $[b,C]\neq\emptyset$, if there is $\bar{C}\in\mathcal{R}_{\mathcal{U}}$ such that $(b,\bar{C})\in\mathcal{D}$, we must have that $(b,C)\in\mathcal{D}$. Let 

$$\mathcal{X} = \{D\in\mathcal{R} : D\leq C \rightarrow \exists b\in\mathcal{AR}(D)\ \ (a\sqsubset b\ \&\ \  (b,C)\in\mathcal{D})\}.$$

Notice that $\mathcal{X}$ is a metric open subset of $\mathcal{R}$. Then by lemma \ref{lemmaOpen} there exists $\bar{C}\in\mathcal{R}_{\mathcal{U}}\cap[depth_C(a),C]$ such that $[a,\bar{C}]$ decides $\mathcal{X}$. Pick $A\in\mathcal{R}_{\mathcal{U}}\cap[a,\bar{C}]$. 

Since $\mathcal{D}$ is dense, there is $(a',A') \leq (a,A)$ such that $(a',A')\in\mathcal{D}$. Notice that $a\sqsubseteq a'$ and hence, as we pointed out in the first parragraph of this proof, we have that $(a',C)\in\mathcal{D}$. And by the definition of $\mathcal{X}$ we also have that $A'\in\mathcal{X}$. But notice that $A'\in[a,A]\subseteq[a,\bar{C}]$; hence we have $[a,\bar{C}]\subseteq\mathcal{X}$ (rather than $[a,\bar{C}]\subseteq\mathcal{X}^{c}$).

Now, it is easy to show that $A$ captures  $(a,\mathcal{D})$: this follows from the fact that $[a,A]\subseteq [a,\bar{C}]\subseteq [a,C]$, and the definition of $\mathcal{X}$.
\end{proof}

\bigskip

\begin{thm}\label{thmGeneric}
Let $\mathcal{U}$ be a selective first-approximation ultrafilter in a given transitive model $M$ of $ZF+DCR$. Then forcing over $M$ with $\mathbb{P}_{\mathcal{U}}$ adds a generic $g\in\mathcal{R}$ with the property that $g\leq^* A$ for every $A\in\mathcal{R}_{\mathcal{U}}$. In fact, every $B\in\mathcal{R}$ is $\mathbb{P}_{\mathcal{U}}$-generic over $M$ if and only if $B\leq^* A$ for every $A\in\mathcal{R}_{\mathcal{U}}$. Also, $M[\mathcal{U}][g] = M[g]$.
\end{thm}
\begin{proof} Suppose $B\in\mathcal{R}$ is $\mathbb{P}_{\mathcal{U}}$-generic over $M$. For every $A\in\mathcal{R}_{\mathcal{U}}$, the set $\{(c,C)\in\mathbb{P}_{\mathcal{U}} : C\leq^*A\}$ is dense open and in $M$: 

Fix $(a,A')\in\mathbb{P}_{\mathcal{U}}$. Take $C_0\in\mathcal{R}_{\mathcal{U}}$ such that $C_0\leq A,A'$. Let $c = r_1(C_0)$. We can assume $depth_{A'}(c)\geq depth_{A'}(a)$ without a loss of generality. Then, by A5, there is $C_1\in [depth_{A'}(c),A']$ such that $\emptyset\neq [c,C_1]\subseteq [c,C_0]$. Notice that $[c,C_1]\subseteq [c,A]$ and hence $C_1\leq^* A$. Also notice that $[a,C_1]\neq\emptyset$ and so $(a,C_1) \leq (a,A')$.

So by genericity, there is one such $(c,C)$ satisfying $B\in [c,C]$. Hence $B\leq^* A$.

\medskip

Now suppose $B\in\mathcal{R}$ is such that $B\leq^* A$ for every $A\in\mathcal{R}_{\mathcal{U}}$. Let $\mathcal{D}$ be a dense open subset of $\mathbb{P}_{\mathcal{U}}$. We need to find $(a,A)\in\mathcal{D}$ such that $B\in [a,A]$. For every $n\in\omega$ find $A_n\in\mathcal{R}_{\mathcal{U}}$, in $M$, such that $A_n$ captures $(r_n(B),\mathcal{D})$. Since $\mathcal{U}$ is in $M$ and selective, we can take $A\in\mathcal{R}_{\mathcal{U}}$, in $M$, such that $A\leq A_n$ for every $n$. By our assumption on $B$ we have that $B\leq^* A$ and so there is $a\in\mathcal{AR}$ such that $[a,B]\subseteq [a,A]$. Let $m=depth_B(a)$. Since $\mathcal{U}$ is closed under finite changes, we can assume without a loss of generality that $a = r_m(B) = r_m(A)$, and so $B\in [m,A]$.  This also implies $A\in [r_m(B),A_m]$ and so $A$ captures $(r_m(B),\mathcal{D})$. Hence, the following is true in $M$: 

\begin{equation}
\forall C\in [m,A]\ \exists n > m\ \ (r_n(C),A)\in\mathcal{D}.
\end{equation}

\medskip

Let $$\mathcal{F} = \{b : \exists n > m\ \  b\in r_n[m,A]\ \&\  (b,A)\not\in\mathcal{D}\},$$ and give to $\mathcal{F}$ the strict end-extension ordering $\sqsubset$. Then the relation $(\mathcal{F},\sqsubset)$ is in $M$, and equation (1) above says that $(\mathcal{F},\sqsubset)$ is well-founded. So by an argument due to Mostowski, (1) holds in the universe and hence, since $B\in [m,A]$, there exists $n > m$\ \ such that $(r_n(B),A)\in\mathcal{D}$. But $B\in [r_n(B),A]$, so $B$ is $\mathbb{P}_{\mathcal{U}}$-generic over $M$.

\medskip

Finally, if $g$ is $\mathbb{P}_{\mathcal{U}}$-generic over $M$ then  $M[\mathcal{U}][g] = M[g]$ follows from the fact that  $\mathcal{R}_{\mathcal{U}}$ (and hence $\mathcal{U}$) can be reconstructed from $g$, in M, as $$\mathcal{R}_{\mathcal{U}} = \{A\in\mathcal{R} : A\in M\ \&\ \ g\leq^* A\}.$$
\end{proof}

\begin{cor}\label{cor2}
If $B$ is $\mathbb{P}_{\mathcal{U}}$-generic over $M$ and $A\leq B$ then $A$ is also $\mathbb{P}_{\mathcal{U}}$-generic over $M$.
\end{cor}
\qed

\medskip

\begin{lemma}\label{lemmaIter}
Let $\breve{\mathcal{U}}$ be the canonical $\mathbb{A}$-name for the $\mathbb{A}$-generic first-approximation ultrafilter. The iteration $\mathbb{A}*\mathbb{P}_{\breve{\mathcal{U}}}$ is equivalent to the abstract Mathias forcing $\mathbb{M}$.
\end{lemma}
\begin{proof} 
Recall that
$$\mathbb{A}*\mathbb{P}_{\breve{\mathcal{U}}} = \{(B,(\dot{a},\dot{A})) : B\in\mathcal{R}\ \&\ \ B\Vdash_{\mathbb{A}} (\dot{a},\dot{A})\in\mathbb{P}_{\breve{\mathcal{U}}}\}.$$ 

Let us see that the mapping $$(a,A)\ \  \rightarrowtail\ \ \ (A,(\breve{a},\breve{A}))$$ is a dense embedding (see \cite{Go}) from $\mathbb{M}$ into $\mathbb{A}*\mathbb{P}_{\breve{\mathcal{U}}}$ (here  $\breve{a}$ and $\breve{A}$ are the canonical $\mathbb{A}$-names for $a$ and $A$, respectively):

\vspace{.25 cm}

 It is easy to show that this mapping preserves the order. 

\vspace{.25 cm}

So given $(B,(\dot{a},\dot{A}))\ \in\ \mathbb{A}*\mathbb{P}_{\breve{\mathcal{U}}}$, we need to find $(d,D)\in\mathbb{M}$ such that $(D,(\breve{d},\breve{D}))\ \leq\ (B,(\dot{a},\dot{A}))$. Since $\mathbb{A}$ is $\sigma$-closed, there exist $a\in\mathcal{AR}$ and $A\in\mathcal{R}$ and also $C\leq^* B$ in $\mathcal{R}$ such that $C\Vdash_{\mathbb{A}} (\breve{a}=\dot{a}\ \ \&\ \ \breve{A}=\dot{A})$\ \ (so we can assume $a\in\mathcal{AR}(C)$). Notice that $(C,(\breve{a},\breve{A}))\ \in\ \mathbb{A}*\mathbb{P}_{\breve{\mathcal{U}}}$ and also that $(C,(\breve{a},\breve{A}))\ \leq\ (B,(\dot{a},\dot{A}))$. Notice as well that  $C\Vdash_{\mathbb{A}}\breve{C}\in\mathcal{R}_{\breve{\mathcal{U}}}$ (since $(C,(\breve{a},\breve{C}))\ \in\ \mathbb{A}*\mathbb{P}_{\breve{\mathcal{U}}}$) and $C\Vdash_{\mathbb{A}}\breve{A}\in\mathcal{R}_{\breve{\mathcal{U}}}$ (since $(C,(\breve{a},\breve{A}))\ \in\ \mathbb{A}*\mathbb{P}_{\breve{\mathcal{U}}}$).

Then, $$C\Vdash\exists x\in\mathcal{R}_{\breve{\mathcal{U}}}\ \ (x\in [\breve{a},\breve{A}]\ \ \&\ \ x\in [\breve{a},\breve{C}]).$$ So there is $D\in\mathcal{R}$ such that $D\in [a, A]\cap [a,C]$. Hence $$(D, (\breve{a},\breve{D}))\ \ \leq\ \  (B,(\dot{a},\dot{A})).$$
\end{proof}

\begin{cor}\label{GenericOmega1Preserve}

Let $M$ be  a given transitive model of $ZF+DCR$, and let $g$ be $\mathbb{M}$-generic over $M$. The following are true:

\begin{enumerate}
\item{Every $A\leq g$ is also $\mathbb{M}$-generic over $M$.}
\item{ $\omega_1^M$ is uncountable in $M[g]$.}
\end{enumerate}

\end{cor}

\begin{proof}

\medskip

1. This follows inmediately from lemma \ref{lemmaIter} and corolary \ref{cor2}.

\medskip

2. By lemma \ref{lemmaIter}, $\mathbb{A}*\mathbb{P}_{\breve{\mathcal{U}}}$ is equivalent to $\mathbb{M}$. So the result follows from the fact that forcing with $\mathbb{A}$ adds no new reals and $\mathbb{P}_{\mathcal{U}}$ satisfies the countable chain condition.
\end{proof}

\section{Solovay models and the abstract Mathias forcing.}

We say that $M$ is a \textit{Solovay model over} $V$ if $M = L(\R)$, where $\R$ is the set of reals in $V^{Coll(\omega,<\kappa)}$, the generic extension of $V$ obtained using the Levy order $Coll(\omega,<\kappa)$ to collapse an inaccesible cardinal $k$ to $\omega_1^{M}$.

In \cite{DiT} there is a proof of the following general result, that we are going to use, about Solovay models:

\begin{prop}\label{prop1}
If $L(\mathbb{R})$ and $L(\mathbb{R}^*)$ are two Solovay models over the same ground model $V$ and $\mathbb{R}\subseteq\mathbb{R}^*$, then $L(\mathbb{R})$ is elementarily embeddable into $L(\mathbb{R}^*)$.
\end{prop}
\qed

\bigskip

Now let $V[G] = V^{Coll(\omega,<\kappa)}$, the Levy collapse of $\kappa$ to $\omega_1$. Let $\R$ be the set of reals in $V[G]$ and let $L(\R)$ be the corresponding Solovay model. Let us force now over $L(\R)$ with $\mathbb{M}$, to add an abstract Mathias real $g$. Let $\R^*$ be the set of reals in $L(\R)[g]$. We are going to show that $L(\R^*)$ is also a Solovay model (see proposition \ref{propExtensionIsSolovay} below). First, we need to prove the following abstract version of a result of Mathias \cite{Mat}:

\begin{defn}
A function $F: \mathcal{X}\subseteq\mathcal{R}\rightarrow [\mathcal{AR}]^{\infty}$ is \textbf{Ramsey-measurable} if  $\forall c\in\mathcal{AR}$ the set $\{D\in\mathcal{X} : c\in F(D)\}$ is Ramsey.
\end{defn}

\begin{lemma}

Given $A\in\mathcal{R}$ and $a\in\mathcal{A}$, let $F:[a,A]\rightarrow [\mathcal{AR}(A)]^{\infty}$ be a Ramsey-measurable function . Then, there exists $B\in [a,A]$ and a function $f:\mathcal{AR}(B)\rightarrow [\mathcal{AR}(A)]^{<\infty}$ such that:

\begin{enumerate}
\item[{(1)}]$c\in f(b) \rightarrow depth_A(c)\leq |b|$.
\item[{(2)}]($\forall b,b'\in \mathcal{AR}(B))\ \ b\sqsubseteq b'\ \rightarrow \ \ \{c\in f(b') : depth_A(c)\leq |b|\} = f(b)$. 
\item[{(3)}]($\forall C\in[a,B]$)($\forall n\in\omega) \{c\in F(C) : depth_A(c)\leq n\} = f(r_n(C))$.

\end{enumerate}

\end{lemma}

\begin{proof}
We will assume $a = \emptyset$ without a loss of generality, to make the proof notationally simpler. So given a Ramsey-measurable function $F:[0,A]\rightarrow [\mathcal{AR}(A)]^{\infty}$, we will find $B\leq A$ and a function $f:\mathcal{AR}(B)\rightarrow [\mathcal{AR}(A)]^{<\infty}$ satisfying $(1)$, $(2)$ and $(3)$ (with $a = \emptyset$).

\medskip

For every $c\in\mathcal{AR}(A)$ consider the following Ramsey set: $$\mathcal{X}_c = \{D\leq A : c\in F(D)\}.$$

\begin{claim}
There exists $B\leq A$ such that for every $b\in\mathcal{AR}(B)$ the following property holds: $$\forall c\in\mathcal{AR}(A)\ \ depth_A(c)\leq |b| \rightarrow [b,B]\ \mbox{decides}\ \mathcal{X}_c.$$
\end{claim}

\begin{proof}[Proof of claim]
We are going to build a fusion sequence $[n,B_n]$ as follows. Let $B_0 = A$.

Suppose $B_n$ has been defined. Let $b_0, b_1, \dots, b_r$ be a list of all the $b$'s in $\mathcal{AR}(B_n)$ with $depth_{B_n} = n$. Find $B^0_n\in [n,B_n]$ such that $[b_0,B^0_n]$ decides $\mathcal{X}_c$, for all $c\in\mathcal{AR}(A)$ with $depth_A(c)\leq |b_0|$.

Suppose $B^k_n\in [n,B_n]$ has been defined, for $0\leq k<r$, such that $[b_k,B^k_n]$ decides  $\mathcal{X}_c$, for all $c\in\mathcal{AR}(A)$ with $depth_A(c)\leq |b_k|$. Notice that $B^k_n$ is compatible with $b_{k+1}$. Then, take  $B^{k+1}_n\in [n,B^k_n]$ such that $[b_{k+1},B^{k+1}_n]$ decides $\mathcal{X}_c$, for all $c\in\mathcal{AR}(A)$ with $depth_A(c)\leq |b_{k+1}|$. 

Let $B_{n+1} = B^r_n$. Then, for every $b\in\mathcal{AR}(B_{n+1})$ with $depth_{B_{n+1}}(b) = n$ we have that $[b,B_{n+1}]$ decides $\mathcal{X}_c$, for all $c\in\mathcal{AR}(A)$ with $depth_A(c)\leq |b|$. Notice that $B_{n+1}\in [n,B_n]$. This concludes the construction of the sequence $[n,B_n]$.

Finally, take $B\in\bigcap_n[n,B_n]$. It is easy to see that for every $b\in\mathcal{AR}(B)$ and every $c\in\mathcal{AR}(A)$ with $depth_A(c)\leq |b|$, $[b,B]$ decides $\mathcal{X}_c$. 
\end{proof}

\bigskip

Let $B\leq A$ be as in the claim and for every $b\in\mathcal{AR}(B)$ choose $B_b\in[b,B]$. Then, define the required function $f:\mathcal{AR}(B)\rightarrow [\mathcal{AR}(A)]^{<\infty}$ as 

$$f(b) = \{c\in F(B_b) : depth_A(c)\leq |b|\}.$$

Now consider any $C\leq B$ and arbitrary $n\in\omega$. Let $b = r_n(C)$, and let $c\in\mathcal{AR}(A)$ be such that $depth_A(c)\leq |b|$. Since $C,B_b\in [b,B]$ and $[b,B]$ decides $\mathcal{X}_c$, we have the following

$$c\in F(C)\ \mbox{iff}\ \ c\in F(B_{b}).$$

This implies that $$\{c\in F(C) : depth_A(c)\leq |b|\} = \{c\in F(B_b) : depth_A(c)\leq |b|\}$$ i.e.,

$$\{c\in F(C) : depth_A(c)\leq n\} = f(r_n(C)).$$

\end{proof}

\begin{cor}\label{corName}
Let $\tau$ be a name for a real number, i.e., a name for an infinite subset of $\mathcal{AR}(A)$ for some $A\in\mathcal{R}$, in the forcing notion $\mathbf{\mathbb{M}}$. Then there exists $B\leq A$ and a function $f:\mathcal{AR}(B)\rightarrow [\mathcal{AR}(A)]^{<\infty}$ such that:

\begin{enumerate}

\item[{(1)}]$\forall b\in\mathcal{AR}(B)\ \ max\{depth_A(c) : c\in f(b)\} \leq |b|$ 
\item[{(2)}]For every $b\in\mathcal{AR}(B)$, $(b,B)$ forces that $f(b)$ is the corresponding initial segment of $\tau$, i.e., $$\{c\in\tau : depth_A(c)\leq |b|\} = f(b).$$ 

\end{enumerate}
\end{cor}
\qed

\bigskip

In virtue of the preceding result and the results of the previous section, the following two lemmas and proposition can be easily proven adapting the proofs of lemmas 2.2 and 2.3 and proposition 2.4 of \cite{DiT}, respectively. In the context of \cite{DiT}, $\mathbb{N}^{[\infty]}$ is the Ramsey space, and $\mathbb{M}$ reduces to the well known (non-abstract) Mathias poset.  

\begin{lemma}[Abstract Local Uniformization Lemma.] Suppose $L(\mathbb{R})$ is a Solovay model. Let $S\subseteq\mathcal{R}\times\mathbb{R}$ be a relation in $L(\mathbb{R})$ such that $\forall x\exists y S(x,y)$. Then, for every nonempty basic open set $[a,A]$ there is $B\in [a,A]$ and a continuous function $h : [a,B] \rightarrow \mathbb{R}$ such that  $\forall x\in [a,B]  S(x,h(x))$.
\end{lemma}
\begin{proof}
We adapt the argument from \cite{DiT}. Take $(a,A)\in\mathbb{M}$ and let $p$ be the real parameter in the definition of $S$. Let $\alpha<\kappa$ be such that $p, [a,A]\in L[G_{\alpha_{0}}]$ (where $G_{\alpha_{0}} = G\cap Coll(\omega,\leq\alpha_{0})$). We may assume $V = L[G_{\alpha_{0}}]$ without a loss of generality. Since $\mathbb{M}$ can be viewed as a subalgebra of $Coll(\omega,<\kappa)$, there is a $V$-generic $g$ for $\mathbb{M}$ in $V^{Coll(\omega,<\kappa)}$. 

Let $y$ be such that $S(g,y)$ hold (in $V^{Coll(\omega,<\kappa)}$). Let $\alpha$ be large enough such that $g,y\in V[G\cap\alpha]$.

Now, take an $\mathbb{M}$-name $\dot{\mathcal{D}}$ for the quotient algebra $Coll(\omega,\alpha)/g$ (in $V[g]$).  Then there is an $\mathbb{M}$-name $\tau$ for a $\dot{\mathcal{D}}$-name for $y$. 

Let $B$ be any $V$-generic real of $\mathbb{M}$ such that $B\in [a,A]$. Define $h:[a,B] \rightarrow \R$ by $h(x)=int_{G_{x}}(\tau)$, where $G_{x}$ is a generic subset of $int_{x}(\dot{\mathcal{D}})$ (in $V[G\cap\alpha]$). By corolary \ref{GenericOmega1Preserve}, each $x\in [a,B]$ is also $\mathbb{M}$-generic over $V$. Note that such $G_{x}$ exists since the collection of all names for subsets of $int_{\dot{\mathcal{D}}}$ can be enumerated in $V[G\cap\alpha]$, and moreover $G_{x}$ can be chosen uniformly in $x$. Therefore, using the real which codes the enumeration of names as a parameter, we have defined a function in $L(\mathbb{R})$ uniformizing the relation $S$: for all $x\in [a,B]$, $S(x,h(x))$ holds because if $x\in [a,B]$ then $[a,B]$ forces $S(\dot{g},\tau)$ and thus, by the Forcing Theorem, $V[G\cap\alpha]\vDash S(x,h(x))$.

Finally, since $\mathcal{R}$ is a topological Ramsey space (i.e., it satisfies the abstract Ellentuck theorem), we can find $B_0\in [a,B]$ such that $f$ is continuous on $[a,B_0]$.
\end{proof}

\begin{lemma}
Let $L(\mathbb{R})$ be a Solovay model. Let $\mathbb{M}$ be the abstract Mathias ordering in  $L(\mathbb{R})$ and $L(\mathbb{R})[g]$ be the corresponding forcing extension. Then, for every formula $\phi$ with parameters from $L(\mathbb{R})$ and quantification only over reals, if $L(\mathbb{R})$ satisfies  $\phi$ so does $L(\mathbb{R})[g]$.
\end{lemma}
\begin{proof}
As usual, the proof is by induction on the complexity of the formula $\phi$. Suppose $L(\mathbb{R})\vDash\phi$ and $L(\mathbb{R})[g]\not\vDash\phi$, where $\phi = \forall x\exists y\psi(x,y)$.

Then, there is a real $x\in L(\mathbb{R})[g]$ such that $L(\mathbb{R})[g]\vDash \forall y\neg\psi(x,y)$. Let $\tau$ be an $\mathbb{M}$-name for $x$ and let $F: [a,A] \rightarrow \R$ be a continuous function associated to $\tau$. Then $(a,A)$ forces $\forall y\neg\psi(F(\dot{g}),y)$.

Let $S = \{(x,y):\psi(F(x),y)\}$. By the abstract local uniformization lemma, there is $B\in [a,A]$ and a continuous function $h : [a,A] \rightarrow \R$ such that $\forall x\in [a,B]\ S(x,h(x))$, i.e., $\forall x\in [a,B]\ \psi(F(x),h(x))$. By the inductive hypothesis, this formula holds in $L(\mathbb{R})[g]$. In particular, $V[G]\vDash\psi(F(\dot{g}),h(\dot{g}))$. So there is a condition in the generic real forcing $\psi(F(\dot{g}),h(\dot{g}))$. This is a contradiction.
\end{proof}

\begin{prop}\label{propExtensionIsSolovay}
If $L(\mathbb{R})$ is a Solovay model over some ground model $V$, and if $g$ is an abstract Mathias real over $L(\mathbb{R})$, then $L(\mathbb{R})[g]$ is also a Solovay model.
\end{prop}
\begin{proof}
In virtue of part 2 of corollary \ref{GenericOmega1Preserve} and lemma 1.1 of \cite{DiT}, we only need to show that every real in $L(\mathbb{R})[g]$ is generic over $V$ by a partial order of size $< \kappa = \omega_1^{L(\mathbb{R})}$. This statement can be expressed by the formula which says that for every two reals $a$ and $x$, there is a real $y$ which codes 

\begin{itemize}
\item[{(a)}]an $\alpha<\omega_1$.
\item[{(b)}]a $Coll(\omega,\alpha)$-name $\tau\in L[a]$.
\item[{(c)}]a $L[a]$-generic subset $G_y$ of $Coll(\omega,\alpha)$ such that  $int_{G_y}(\tau) = x$.
\end{itemize}

Since this is a formula with quantification over the reals which has only reals as parameters, by the previous lemma it is absolute between $L(\mathbb{R})$ and  $L(\mathbb{R})[g]$. Since  $L(\mathbb{R})$ is a Solovay model, it satisfies this formula, and therefore the $\mathbb{M}$-extension satisfies it as well, and hence it is a Solovay model.	
\end{proof}

\bigskip

So by proposition \ref{prop1}, we have the following:

\begin{cor}\label{corEmbed}
If $L(\mathbb{R})$ is a Solovay model over some ground model $V$, and if $g$ is an abstract Mathias real over $L(\mathbb{R})$, then $L(\mathbb{R})$ is elementarily embeddable into  $L(\mathbb{R})[g]$.
\end{cor}
\qed

\bigskip

\section{The Open Coloring Axiom in $L(\mathbb{R})[\mathcal{U}]$.}\label{sect6}

In this section we show how to extend the result of \cite{DiT} about selective ultrafilters on $\N$ to the context of first-approximation ultrafilters.

\bigskip

For a set $X$ of reals, let $OCA(X)$ be the following statement:

\medskip

Given an open subset $K$ of $\mathbb{R}^{[2]}$,  either there is a perfect $P\subseteq X$ such that $P^{[2]}\subseteq K$, or there is a sequence $\{X_i : i\in\omega\}$ such that $X = \bigcup_i X_i$ and $X_i^{[2]}\cap K = \emptyset$, for every $i$. If the first alternative holds we say that $X$ contains a $K$-\textit{perfect} set. If the second alternative holds we say that $X$ is $K$-\textit{countable}.

\medskip

Recall that \textit{Open Coloring Axiom}, $OCA$, states that $OCA(X)$ holds for every set $X$ of real numbers (see \cite{Tod2}). This definitions are naturally extended to separable metric spaces in general.

\medskip

Different  proofs of the following important result that we are going to use can be found in \cite{DiT} and \cite{Fen}:

\begin{thm}
If  $L(\mathbb{R})$ is a Solovay model then OCA($X$) holds for every set of reals $X$ in $L(\mathbb{R})$.
\end{thm}
\qed

\bigskip

The following is our abstract version of lemma 5.2 of \cite{DiT}, which turns out to be an impotant tool to the proof of OCA in $L(\mathbb{R})[\mathcal{U}]$.

\medskip

\begin{lemma}
Let $M$ be a given ground model and let $\mathcal{U}$ be a selective first-approximation ultrafilter in $M$. Let $K$ be, in $M$, a subset of $\mathbb{R}^{[2]}$. If a set of reals $X$ in $M$ is $K$-countable in the  $\mathbb{P}_{\mathcal{U}}$-generic extension then $X$ is $K$-countable in $M$. 
\end{lemma}
\begin{proof} Let $\breve{K}$ and  $\breve{X}$ be the canonical names for $K$ and $X$ respectively, and suppose there is $(a,A)\in\mathbb{P}_{\mathcal{U}}$ forcing that $\breve{X}$ is $\breve{K}$-countable. Let $(\dot{X}_n)_n$ be a sequence of names such that $(a,A)$ forces $\breve{X} = \bigcup_n\dot{X}_n$ and $\dot{X}_n^{[2]}\cap\breve{K} = \emptyset$ for all $n\in\omega$.

For an integer $n$ and $b\in\mathcal{AR}$ with $a\sqsubseteq b$, let $$X_{(n,b)} = \{x\in X : (\exists B\in\mathcal{U})[(b,B)\Vdash x\in\dot{X}_n\}.$$ Notice that $$X = \bigcup\{X_{(n,b)} : n\in\omega, b\in\mathcal{AR}, a\sqsubseteq b\}.$$ Let us see that $X_{(n,b)}\cap K = \emptyset$ for all $n$ and $b$: consider $x_0,x_1\in X_{(n,b)}$ such that $x_0\neq x_1$. Let $B_0,B_1\in\mathcal{R}_{\mathcal{U}}$ be such that $(b,B_i)\Vdash x_i\in\dot{X}_n$ ($i\in\{0,1\}$). Since $\mathcal{U}$ is a first-approximation ultrafilter, there exist $C\in\mathcal{R}_{\mathcal{U}}$ such that $C\in [b,B_0]\cap [b,B_1]\cap [b,A]$. Then $(b,C)\Vdash (x_0,x_1)\in X_{(n,b)}^{[2]}\ \ \mbox{and also}\ \ (b,C)\Vdash (x_0,x_1)\not\in K$.
\end{proof}

\begin{thm}
Let $L(\mathbb{R})$ be a Solovay model and let $\mathcal{U}$ be a selective first-approximation ultrafilter generic over $L(\mathbb{R})$. Then $L(\mathbb{R})[\mathcal{U}]$ satisfies OCA.
\end{thm}
\begin{proof} 

As in section \ref{sec4}, we force over $L(\mathbb{R})$ with $\mathbb{A}$ to add a selective first-approximation ultrafilter. In the extension, consider the poset $\mathbf{\mathbb{P}_{\mathcal{U}}}$ and let $g$ be $\mathbb{P}_{\mathcal{U}}$-generic over $L(\mathbb{R})[\mathcal{U}]$. Then $L(\mathbb{R})[\mathcal{U}][g] = L(\mathbb{R})[g]$ by theorem \ref{thmGeneric}. Also recall that the iteration $\mathbb{A}*\mathbb{P}_{\mathcal{U}}$ is equivalent to the abstract Mathias forcing $\mathbf{\mathbb{M}}$ (see lemma \ref{lemmaIter}). Hence the real $g$ is an abstract Mathias real over $L(\mathbb{R})$. 

Let $\mathbb{R}^*$ be the set of reals in the extension $L(\mathbb{R})[g]$, and consider $L(\mathbb{R}^*) = L(\mathbb{R})[g]$. By corollary \ref{corEmbed}, there exists an elementary embedding $j:L(\mathbb{R})\rightarrow L(\mathbb{R}^*)$ which fixes the reals of $L(\mathbb{R})$ and the ordinals.

Let $X,K\in L(\mathbb{R})[\mathcal{U}]$ be such that, in $L(\mathbb{R})[\mathcal{U}]$, $X$ is a set of reals and $K$ is an open subset of $\mathbb{R}^{[2]}$. We will use the embedding $j$ to prove OCA$(X)$ in $ L(\mathbb{R})[\mathcal{U}]$.

Let $\dot{X}$ be a name for $X$. Then $j(\dot{X})$ is a name for a set of reals in  $L(\mathbb{R}^*)$ with the same definition. Since  $L(\mathbb{R})$ and  $L(\mathbb{R})[\mathcal{U}]$ have the same reals we have that $K$ is in  $L(\mathbb{R})$. This implies that $j(K)$ is open and also, it is coded by the same real as $K$.

Assume, without a loss of generality, that every condition forces that $\dot{X}$ is not $K$-countable. Consider $$Y = \{y\in\mathbb{R}^* : g\Vdash y\in j(\dot{X})\},$$ which is a set in $L(\mathbb{R}^*)$. Since $L(\mathbb{R}^*)$ is a Solovay model, applying OCA (in $L(\mathbb{R}^*)$) to $Y$ and $j(K)$ we have that either $Y$ is  $j(K)$-countable in $L(\mathbb{R}^*)$, or there is a perfect set $P\subseteq Y$ such that $P^{[2]}\subseteq j(K)$.

Notice that the interpretation of $X$ in  $L(\mathbb{R}^*)$ is a subset of $Y$: if $x\in X$ then there is $A\in\mathcal{R}_{\mathcal{U}}$ such that, in $L(\mathbb{R})$, $A\Vdash x\in\dot{X}$. Then by elementarity, in  $L(\mathbb{R}^*)$, $A\Vdash x\in j(\dot{X})$. But $g\leq^* A$  by theorem \ref{thmGeneric} and so $g\Vdash x\in j(\dot{X})$.

Then, by our assumption on $X$ and by the previous lemma,  $Y$ cannot be $j(K)$-countable in $L(\mathbb{R}^*)$. So let $P$ be a perfect set witnessing OCA($Y$). By the definition of $Y$, we have $g\Vdash P\subseteq j(\dot{X})$. Therefore, 

\begin{flushleft}
$$L(\mathbb{R}^*)\vDash(\exists A\in\mathcal{R})(\exists P\in\mathbb{R})[P\ \ \mbox{is a}\ \ j(K)\mbox{-perfect set and}\ \ A\Vdash P\subseteq j(\dot{X})]$$
\end{flushleft}

 and then, 
\begin{flushleft}
$$L(\mathbb{R})\vDash(\exists A\in\mathcal{R})(\exists P\in\mathbb{R})[P\ \ \mbox{is a}\ \ K\mbox{-perfect set and}\ \ A\Vdash P\subseteq \dot{X}]$$
\end{flushleft}

 by elementarity. So $X$ is $K$-perfect in $L(\mathbb{R})[\mathcal{U}]$.

\end{proof}

\begin{cor}
If $L(\mathbb{R})$ be a Solovay model and  $\mathcal{U}$ is a selective first-approximation ultrafilter generic over $L(\mathbb{R})$ then $L(\mathbb{R})[\mathcal{U}]$ satisfies the perfect set property.
\end{cor}
\begin{proof}
Setting $K = \mathbb{R}^{[2]}$, $OCA(X)$ reduces to ``$X$ is countable or it contains a perfect set''.
\end{proof}

\medskip

One interesting final comment on the results of this section is that they are really ``scheme-results'', in the sense that we have provided simultaneous proofs for the corresponding results in each particular topological Ramsey space by means of our abstract definition of the ultrafilter $\mathcal{U}$.

\bigskip

\begin{ack}
The author would like to show his gratitude to the following institutions and persons:  Consejo de Desarrollo Cient\'ifico y Human\'istico (CDCH), Coodinaci\'on de Investigaci\'on de la Facultad de Ciencias, Coodinaci\'on de Postgrado de Matem\'aticas and Escuela de Matem\'aticas of the Universidad Central de Venezuela (UCV), for finantial and institutional support. Stevo Todorcevic for financial support through his grant of the National Science Foundation of Canada (NSFC) and also for suggesting the author to venture in this research and for supervising its development in all of its stages. Ida Bullat and the staff of the Department of Mathematics of the University of Toronto. Ilijas Farah and the Set theory Seminar at Fields Institute. Carlos Di Prisco and the Seminario de Logica Matem\'atica IVIC-UCV for many years of invaluable teachings.
\end{ack}

\end{document}